\documentclass[12pt]{article}
\usepackage{amssymb}
\usepackage{latexsym}

\oddsidemargin=\evensidemargin
\begin{document}
\newfont{\blb}{msbm10 scaled\magstep1} 

\newtheorem{defi}{Definition}
\newtheorem{theo}{Theorem}[section]
\newtheorem{theo1}{Theorem}[section]
\newtheorem{prop}[theo]{Proposition}
\newtheorem{lemm}[theo]{Lemma}
\newtheorem{coro}[theo]{Corollary}
%\markright{Torsion free groups with all subgroups $4$-subnormal}
\date{}
\author{Gunnar Traustason \\
James Williams \\
Department of Mathematical Sciences \\ 
University of Bath, UK}
\title{Powerfully nilpotent groups}
\maketitle
\begin{abstract}
\mbox{}\\
We introduce a special class of powerful $p$-groups that  we call {\it powerfully
nilpotent groups} that are finite $p$-groups that possess a central series of 
a special kind. To these we can attach the notion of a powerful nilpotence class that leads naturally to a classification in terms of  an `ancestry tree' and powerful coclass. We show that there are finitely many powerfully nilpotent $p$-groups of each given powerful coclass and develop some general theory for this class of groups. We also determine the growth of powerfully nilpotent groups of exponent $p^{2}$ and order $p^{n}$ where $p$ is odd. The number of these is $f(n)=p^{\alpha n^{3}+o(n^{3})}$ where $\alpha=\frac{9+4\sqrt{2}}{394}$. For the larger class of all powerful groups of exponent $p^{2}$ and order $p^{n}$, where $p$ is odd, the number is $p^{\frac{2}{27}n^{3}+o(n^{3})}$. Thus here the class of powerfully nilpotent $p$-groups is large while sparse within the larger class of powerful $p$-groups. 
\end{abstract}
\section{Introduction and basic properties}
In this paper we study a special subclass of powerful $p$-groups that we call powerfully nilpotent $p$-groups. The notion of a powerfully nilpotent group originates from the theory of  symplectic alternating algebras [4,6]. The latter is a certain class of alternating algebras that were introduced to study powerful $2$-Engel groups [2,5]. There turned out to be a 1-1 correspondence between symplectic alternating algebras over $\mbox{GF}(3)$ and a rich class of powerful $2$-Engel $3$-groups. Under this correspondence, the nilpotent symplectic
alternating algebras correspond to powerful $3$-groups with some special extra property that means having a central series of a certain kind. Singling out this property gives rise to powerfully
nilpotent $p$-groups that turns out to be a rich subclass of powerful $p$-groups that exist for any prime $p$. As we will see the property is natural and leads to a  class of $p$-groups with some beautiful properties. As another motivation we will see that these groups occur frequently as characteristic subgroups of any powerful $p$-group $G$. For example the Frattini
subgroup of $G$ will always be powerfully nilpotent and the same is true for all the proper terms of the derived series and the lower central series of $G$.  Although we do not make any explicit use of GAP [7] in this paper, a number of GAP calculations  provided insight for some of the results.\\ \\ 
\normalsize
Let $G$ be a finite $p$-group where $p$ is a prime. \\ \\
{\bf Definition}. Let $H\leq K\leq G$. An ascending chain of subgroups 
  $$H=H_{0}\leq H_{1}\leq \cdots \leq H_{n}=K$$
is {\it powerfully central} if $[H_{i},G]\leq H_{i-1}^{p}$ 
for $i=1,\ldots ,n$. Here $n$ is called the {\it length} of the chain. \\ \\
{\bf Remark}. Recall that a finite $p$-group is said to be {\it powerful} if either $p$ is odd and $[G,G]\leq G^{p}$ or $p=2$ and $[G,G]\leq G^{4}$. More generally, a subgroup $H$ of $G$ is {\it powerfully embedded} in $G$ if $[H,G]\leq H^{p}$, when $p$ is odd, and $[H,G]\leq H^{4}$ for $p=2$. \\ \\
{\bf Definition}. A powerful $p$-group 
 $G$ is {\it powerfully nilpotent} if it has a powerfully central ascending chain of subgroups of the form
    $$\{1\}=H_{0}\leq H_{1}\leq \cdots \leq H_{n}=G.$$ 
{\bf Remark}. Let $G$ be a powerful $2$-group of exponent $2^{e}$ and consider the ascending chain
    $$\{1\}=G^{2^{e}}\leq G^{2^{e-1}}\leq \cdots \leq G^{2}\leq G.$$
According to Shalev's interchange lemma [3] we know that if $M,N$ are powerfully embedded subgroups of a finite $p$-group $G$, then $[M^{p^{i}},N^{p^{j}}]=[M,N]^{p^{i+j}}$ for all non-negative integers $i,j$. In particular 
$[G^{2^{k}},G]=[G,G]^{2^{k}}\leq G^{2^{k+2}}=(G^{2^{k+1}})^{2}$ and thus the chain above is
powerfully central and therefore $G$ is powerfully nilpotent. Notice also that if $e\geq 2$ then $[G^{2^{e-2}},G]=[G,G]^{2^{e-2}}\leq G^{2^{e}}=\{1\}$ and thus $G^{2^{e-2}}\leq Z(G)$ and the powerful class of $G$ is at most $e-1$. Having seen that all powerful $2$-groups are powerfully nilpotent we will normally assume in the following that a given prime $p$ is odd except if otherwise stated.  \\ \\
{\bf Remark}. For any prime $p$ one sees from the same calculation that if $[G,G]\leq G^{p^{2}}$, then $G$ is powerfully nilpotent. \\ \\
{\bf Definition}. We say that a finite $p$-group $G$ is {\it strongly powerful} if $[G,G]\leq G^{p^{2}}$. \\ \\
{\bf Remarks}. (1) It follows from the discussion above that a strongly powerful
group must be powerfully nilpotent. \\ \\
(2) Let $G$ be a finite $p$-group and let $N$ be a subgroup of the Frattini subgroup $\Phi(G)$ where $N$ is powerfully embedded in $G$. Then
$N$ is strongly powerful and thus powerfully nilpotent. To see this notice that
      $$[N,N]\leq [N,[G,G]G^{p}]\leq [N,G,G][N,G]^{p}\leq [N^{p},G][N,G]^{p}\leq [N,G]^{p}=N^{p^{2}}.$$
This implies in particular that if $G$ is powerful then $\Phi(G)$ and more generally $G^{p^{i}}$ for $i\geq 1$ are powerfully
nilpotent. The same is true for 
the proper terms of the derived and the lower central series.  \\ \\ 
(3)  Let $p$ be an odd prime.  Consider any powerfully nilpotent group with a powerfully
central chain $\{1\}=H_{0}\leq H_{1}\leq \cdots \leq H_{n}=G$. Then in particular
$[H_{i},G]\leq H_{i}^{p}$ and each term of the chain is powerfully embedded in 
$G$. Notice also that each term is powerfully nilpotent. Clearly a quotient
of a powerfully nilpotent group is also powerfully nilpotent. As for powerful
groups the property of being powerfully nilpotent is not in general inherited by subgroups. \\ \\
{\bf Definition}. If $G$ is powerfully nilpotent then the {\it powerful nilpotence class} of $G$ is the shortest length that a powerfully central chain of $G$ can have. \\ \\
{\bf The upper powerfully central series}. This is defined recursively as 
follows: $\hat{Z}_{0}(G)=\{1\}$ and for $n\geq 1$, 
      $$\hat{Z}_{n}(G)=\{a\in G:\, [a,x]\in \hat{Z}_{n-1}(G)^{p}\mbox{ for
all }x\in G\}.$$
Notice in particular that $\hat{Z}_{1}(G)=Z(G)$. \\ \\
{\bf Remark}. The upper powerfully central series is clearly
the fastest ascending series that is powerfully central. There is however no
canonical fastest lowest descending series that is powerfully central. It is easy to see that if $G$ is powerfully nilpotent, then
the powerful nilpotence class is the smallest nonzero integer $n$ where
$\hat{Z}_{n}(G)=G$. \\ \\
If $G/G^{p^{2}}$ is abelian then $G$ is strongly powerful and thus powerfully nilpotent. It turns out that in order for $G$ to powerfully nilpotent, it suffices
that $G/G^{p^{2}}$ is powerfully nilpotent. \\
\begin{prop} Let $G$ be any finite $p$-group of exponent $p^{e}$ where $e\geq 2$. If $G/G^{p^{2}}$ is powerfully
nilpotent, then $G$ is powerfully nilpotent. Furthermore if $G/G^{p^{2}}$ has 
powerful class $m$, then the powerful class of $G$ is at most $(e-1)m$. 
\end{prop}
{\bf Proof}\ \ Notice first that, as $G/G^{p^{2}}$ is powerful, we have $[G,G]\leq G^{p}$ when $p$ is odd and $[G,G]\leq G^{4}$ when $p=2$.  Thus $G$ is powerful. We have seen above that a powerful $2$-group is powerfully nilpotent and that the powerful class is then at most $e-1$. \\ \\
We can thus assume that $p$ is odd.  Suppose 
    $$G/G^{p^{2}}=\overline{H}_{0}\geq \overline{H}_{1}\geq\cdots \geq 
\overline{H}_{m-1}\geq \{\overline{1}\}$$
is the upper powerfully central series for $G/G^{p^{2}}$, where $\overline{H_{i}}=
H_{i}/G^{p^{2}}$ with $G^{p^{2}}\leq H_{i}\leq G$. Since $G$ is powerful, we have $[G^{p},G]\leq G^{p^{2}}$ and thus $G^{p}\leq H_{m-1}$. 
Hence 
  $$G=H_{0}\geq H_{1}\geq\cdots \geq H_{m-1}\geq  H_{m}=G^{p}$$
is powerfully central. Consider the descending chain 
$$\begin{array}{lllllllllll}
\mbox{}    G  & = & H_{0} & \geq & H_{1}  & \geq & \cdots & \geq & H_{m} & = & G^{p} \\
\mbox{}   G^{p} & = & H_{0}^{p} & \geq & H_{1}^{p} & \geq & \cdots & \geq & H_{m}^{p} & = &
   G^{p^{2}} \\ 
\mbox{}     &  &  &  &  & \vdots  &  &  &  &  & \\
\mbox{}G^{p^{e-2}} & = & H_{0}^{p^{e-2}} & \geq & H_{1}^{p^{e-2}} & \geq &
\cdots & \geq & H_{m}^{p^{e-2}} & = & G^{p^{e-1}}.
\end{array}$$
We know already that the first line gives us a powerfully central chain and as a result $H_{0},\ldots ,H_{m}$ are powerfully embedded in $G$. Using this fact and
the interchange lemma, we have that 
$$[H_{i}^{p^{k}},G]=[H_{i},G]^{p^{k}}\leq H_{i+1}^{p^{k+1}}=
(H_{i+1}^{p^{k}})^{p}$$
for $0\leq i\leq m-1$ and $0\leq k\leq e-1$. It follows that we have a powerfully central chain. Notice also that $[H_{m-1}^{p^{e-2}},G]=[H_{m-1},G]^{p^{e-2}}
\leq (G^{p^{2}})^{p^{e-2}}=\{1\}$ and thus $H_{m-1}^{p^{e-2}}\leq Z(G)$. It follows from this and the powerfully central chain above that the powerful
class of $G$ is at most $(e-1)m$. $\Box$
\section{The ancestry tree and powerful coclass}
{\bf Definition}. Let $G$ be a powerfully nilpotent $p$-group of powerful
class $c$ and order $p^{n}$. We define the {\it powerful coclass} of $G$ to be the number $n-c$. \\ \\
A natural approach is to develop something that corresponds to a coclass theory for finite $p$-groups where coclass is replaced by powerful coclass. We are indeed 
going to see that there are only finitely many powerfully nilpotent  
$p$-groups of any given powerful coclass and this will follow from sharp bounds for
the rank and the exponent in terms of the powerful coclass. More precisely, if $G$ is
a powerfully nilpotent $p$-group of rank $r$ and exponent $p^{e}$ then 
$r\leq n-c+1$ and $e\leq n-c+1$. In this section we will prove the first inequality. For the `classification' we introduce, for a given prime $p$, the {\it ancestry tree} where every
powerfully nilpotent $p$-group occurs exactly once. \\ \\
{\bf The ancestry tree}. Let $p$ be a fixed prime. The {\it vertices} of the 
ancestry tree are all the powerfully nilpotent $p$-groups (one for each isomorphism class). Two vertices $G$ and $H$ are joined by a directed edge from $H$ to
$G$ if and only if $H\cong G/Z(G)^{p}$ and $G$ is not abelian. Notice that this implies that $Z(G)^{p}\not =\{1\}$ and thus the powerful class of $G$ is one more than that of $H$. 
We then also say that 
$G$ is a {\it direct ancestor} of $H$ or that $H$ is a {\it direct descendant}
of $G$, and we write $H\rightarrow G$. Notice that a direct descendant is unique. A group $H$ can however have
a number of (even infinitely many) direct ancestors.  \\ \\ 
More generally, if for powerfully nilpotent groups $H$ and $G$ there is a chain
     $$H=H_{0}\rightarrow H_{1}\rightarrow \cdots \rightarrow H_{n}=G,$$
then we say that $G$ is an {\it ancestor} of $H$ and that $H$ is a {\it descendant} of $G$. Notice that if $H$ has powerful class $c$ then $H_{i}$ has powerful
class $c+i$. \\ \\
{\bf Example 1}. For $n\geq 7$ let $G(n):=\langle a,b:\,a^{3^{n}}=1,\, b^{27}=1,\,[a,b]=a^{3^{n-3}}\rangle$. Then $G(n)$ is a direct ancestor to $C^{3^{4}}\times C^{3^{3}}$.
A powerfully nilpotent group can thus have infinitely many direct ancestors. \\ \\
{\bf Example 2}.   For $n\geq 2$ let $G(n):=\langle a,b:\,a^{p^{n}}=1,\,b^{p^{n}}=1\, [b,a]=a^{p^{2}}\rangle$. Notice that $G(n+1)/Z(G)^{p}\cong G(n)$. Thus the
group $C_{p^{2}}\times C_{p^{2}}$ has infinitely many generations of ancestors. \\ \\
{\bf Remark}. Suppose $H$ is a powerfully nilpotent $p$-group of order 
$p^{n(H)}$ and powerful class $c(H)$ and suppose $G$ is a direct ancestor of $H$. Then
$G$ has powerful class $c(G)=c(H)+1$ and order $p^{n(G)}=|G/Z(G)^{p}|\cdot |Z(G)^{p}|=
|H|\cdot |Z(G)^{p}|=p^{n(H)+k}$ where $|Z(G)^{p}|=p^{k}$. Notice that $k\geq 1$. Thus the powerful coclass of $G$ is $d(G)=n(H)+k-(c(H)+1)=(n(H)-c(H))+(k-1)$ and thus 
        $$d(G)\geq d(H)$$
where $d(H)$ is the powerful coclass of $H$. Notice that we have 
equality if and only if $|Z(G)^{p}|=p$. \\ \\
We now turn to the task of proving that the rank of $G$ is bounded in terms
of the powerful coclass. 
\begin{lemm} Let $\{1\}=H_{0}\leq H_{1}\leq\cdots \leq H_{n}=G$ be any
powerfully central chain for a powerfully nilpotent $p$-group $G$ where $|H_{k}|=p^{k}$ for
$k=0,1,\ldots ,n$.  Suppose $H_{j}^{p}=H_{j-1}^{p}$ for some $1\leq j\leq n$. Then there exists
$x\in H_{j}\setminus H_{j-1}$ such that $x^{p}=1$. 
\end{lemm}
{\bf Proof}\ \ If $j=1$ this is obvious. Suppose then that $j\geq 2$.
Pick first any $x\in H_{j}\setminus H_{j-1}$. Then as $H_{j}^{p}=
H_{j-1}^{p}$ we have that $x^{p}=y^{p}$ for some $y\in H_{j-1}$. Now $x$ commutes
with $y$ modulo $H_{j-2}^{p}$. Hence 
     $$(xy^{-1})^{p}=x^{p}y^{-p}z^{p}=z^{p}$$
for some $z\in H_{j-2}$ and by replacing $x$ by $xy^{-1}$ we can now assume that
$x^{p}\in H_{j-2}^{p}$. If $j=2$ we are done, otherwise a similar argument shows that we can replace $x$ by
a new element in $H_{j}\setminus H_{j-1}$ such that $x^{p}\in H_{j-3}^{p}$. Continuing like this we see that we can eventually pick our $x$ such that $x^{p}=1$.
$\Box$ \\  
\begin{prop} Let $\{1\}=H_{0}\leq H_{1}\leq\cdots \leq H_{n}=G$ be any powerfully central chain for a powerfully nilpotent $p$-group $G$ where $|H_{k}|=p^{k}$ for $k=0,\ldots ,n$. We can then choose $a_{1},\ldots ,
a_{n}\in G$ such that $H_{i}=\langle a_{1},\ldots, a_{i}\rangle$ and such that
$a_{i}^{p}=1$ if and only if $H_{i}^{p}=H_{i-1}^{p}$ for $i=1,\ldots ,n$. 
\end{prop}
{\bf Proof}\ \ Follows directly from Lemma 2.1. $\Box$ \\ \\
{\bf Definition}. Let $G$ be a powerfully nilpotent group of order $p^{n}$ and
let ${\mathcal L}$ be an ascending powerfully central chain
   $$\{1\}=H_{0}\leq H_{1}\leq \cdots \leq H_{n}=G$$
where $|H_{j}|=p^{j}$ for $0\leq j\leq n$. The number 
      $$s_{\mathcal L}(G)=|\{H_{0}^{p},H_{1}^{p},\ldots ,
H_{n}^{p}\}|$$
is called the {\it $p$th power length} of ${\mathcal L}$. \\ \\
{\bf Remark}. We can thus think of $s_{\mathcal L}(G)-1$ as the number
of `jumps' in the chain 
     $$\{1\}=H_{0}^{p}\leq H_{1}^{p}\leq \cdots \leq H_{n}^{p}=G^{p}.$$ 
%
%{\bf Remark}. It follows from Lemma 2.1 that we can choose elements 
%$a_{1}\in H_{1}\setminus H_{0},\ a_{2}\in H_{2}\setminus H_{1},\ \ldots
%,a_{n}\in H_{n}\setminus H_{n-1}$ such that 
%
%$$\begin{array}{l}
%%
%    a_{i}^{p}=1\mbox{ if }H_{i}^{p}=H_{i-1}^{p} \\
%
%   a_{i}^{p}\not\in H_{i-1}^{p}\mbox{ if }H_{i}^{p}\not =H_{i-1}^{p}.
%
%\end{array}$$
%
%Notice that the number of nontrivial powers is $s_{\mathcal L}(G)-1$ and that
%$H_{i}=\langle a_{1},\ldots ,a_{i}\rangle$ with $H_{i}^{p}=\langle a_{i}^{p}\ra%ngle H_{i-1}^{p}$ for $i=1,\ldots ,n$. 
%
\begin{lemm} If ${\mathcal L}$ is any ascending powerfully central
chain as above then $|G^{p}|=p^{s_{\mathcal L}(G)-1}$. 
\end{lemm}
{\bf Proof}\ \ Suppose that ${\mathcal L}$ is the asending central 
chain
      $$\{1\}=H_{0}\leq H_{1}\leq \cdots \leq H_{n}=G.$$
Let $t=s_{\mathcal L}(G)-1$ and suppose 
    $$\{1\}=H_{0}^{p}=\ldots =H_{j_{1}-1}^{p}<H_{j_{1}}^{p}=
\ldots =H_{j_{2}-1}^{p}<H_{j_{2}}^{p}\cdots H_{j_{t}-1}^{p}<H_{j_{t}}^{p}=
\cdots =H_{n}^{p}=G^{p}.$$
From Proposition 2.2  we know that we can pick elements $a_{1},\ldots
,a_{n}\in G$ such that $H_{i}=\langle a_{1},\ldots ,a_{i}\rangle$ and
that $a_{j_{k}}^{p}\not\in H_{j_{k}-1}^{p}$. Notice that $a_{j_{k}}^{p}\in H_{j_{k-1}}$
and thus $a_{j_{k}}^{p^{2}}\in H_{j_{k-1}}^{p}$ and thus $H_{j_{k}}^{p}/H_{j_{k-1}}^{p}$ a cyclic group of order $p$. Thus 
    $$|G^{p}|=\frac{|H_{n}^{p}|}{|H_{n-1}^{p}|}\cdots \frac{|H_{1}^{p}|}{|H_{0}^{p}|}=\frac{|H_{j_{t}}^{p}|}{|H_{j_{t}-1}^{p}|}\cdots \frac{|H_{j_{1}}^{p}|}{|H_{j_{1}-1}^{p}|}=p^{t}.$$
Hence the result follows. $\Box$ \\
\begin{prop} Let $G$ be a powerfully nilpotent $p$-group of rank $r$ and
order $p^{n}$. Then
$s_{\mathcal L}(G)=n-r+1$.
\end{prop}
{\bf Proof}\ \ As $G^{p}$ is the Frattini subgroup of $G$ it follows that
      $$p^{r}=|G/G^{p}|=\frac{p^{n}}{p^{s_{\mathcal L}(G)-1}}.$$
The result follows from this. $\Box$ \\ \\
{\bf Remark}\ \ In particular it follows that $s_{\mathcal L}(G)$ does not
depend on ${\mathcal L}$. We will thus denote this number by $s(G)$ and will call it the {\it $p$th power length} of $G$. \\
\begin{lemm} Let $G$ be a powerfully nilpotent group of powerful class $c\geq 2$ then
   $$[G,G]=[\hat{Z}_{c}(G),G]>[\hat{Z}_{c-1}(G),G]>\ldots >[\hat{Z}_{1}(G),G]=\{1\}$$
and
   $$G^{p}=\hat{Z}_{c}(G)^{p}\geq \hat{Z}_{c-1}(G)^{p}>\ldots >\hat{Z}_{1}(G)^{p}>\hat{Z}_{0}(G)^{p}=\{1\}.$$
In particular $|G^{p}|\geq |[G,G]|\geq p^{c-1}$. 
\end{lemm}
{\bf Proof}\ \ Suppose $2\leq j\leq c$. If $[\hat{Z}_{j}(G),G]=[\hat{Z}_{j-1}(G),G]$, then $[\hat{Z}_{j}(G),G]\leq \hat{Z}_{j-2}(G)^{p}$ and thus we get the contradiction that $\hat{Z}_{j}(G)\leq \hat{Z}_{j-1}(G)$. The proof of the latter strict inequalities is similar. Let $1\leq j\leq c-1$. If $\hat{Z}_{j}(G)^{p}=\hat{Z}_{j-1}(G)^{p}$ then $[\hat{Z}_{j+1}(G),G]\leq \hat{Z}_{j-1}(G)^{p}$ and thus $\hat{Z}_{j+1}(G)=\hat{Z}_{j}(G)$ that gives the contradiction that the powerful class of $G$ is at most $j\leq c-1$. $\Box$ 
\begin{theo} Let $G$ be a powerfully nilpotent $p$-group with rank $r$, powerful class $c$ and order $p^{n}$. Then 
        $$r\leq n-c+1.$$
\end{theo}
{\bf Proof}\ \ Refining the upper powerfully central series we get an ascending powerfully central series 
   $$\{1\}=H_{0}<H_{1}<\ldots <H_{n}=G$$
where $|H_{i}|=p^{i}$ for $i=0,\ldots ,n$ and where $\{H_{0}^{p}, \ldots ,H_{n}^{p}\}$ contains $\hat{Z}_{0}(G)^{p}, \hat{Z}_{1}(G)^{p},\ldots ,\hat{Z}_{c-1}(G)^{p}$. By Lemma 2.5 these $c$ groups are distinct and thus $c\leq s(G)=n-r+1$. Hence $r\leq n-c+1$. $\Box$  \\ \\ 
{\bf Presentations of powerfully nilpotent groups}. Let $G$ be any powerfully nilpotent $p$-group of exponent $p^{e}$, order
$p^{n}$ and rank $r$. Suppose that 
    $$G=H_{0}> H_{1}>\cdots > 
H_{m}=G^{p}$$
is a powerfully central series as given in the proof of Proposition 1.1. We can always refine such a powerfully central chain to get a chain of length $r$ such that the factors are of order $p$. Without loss of generality we can thus assume that we
have a powerfully central chain
\begin{equation} 
  G=H_{0}>H_{1}> \cdots > H_{r}=G^{p}
\end{equation}
such that $|H_{i}/G^{p}|=p^{r-i}$ for $i=0,\ldots ,r$. We then have 
$G=\langle a_{1},\ldots ,a_{r}\rangle$ for any choice of elements $a_{1},
\ldots ,a_{r}\in G$ where $H_{i}=\langle a_{i+1},\ldots, a_{r}\rangle G^{p}$ for $i=0,\ldots ,r$. We would like to choose $a_{1},\ldots ,a_{r}$
so that the choice best reflects the structure of the group $G$. \\ \\
As we saw in the proof of Proposition 1.1 we then get a powerfully central series 
$$\begin{array}{lllllllllll}
\mbox{}    G  & = & H_{0} & > & H_{1}  & > & \cdots & > & H_{r} & = & G^{p} \\
\mbox{}   G^{p} & = & H_{0}^{p} & \geq & H_{1}^{p} & \geq & \cdots & \geq & H_{r}^{p} & = &
   G^{p^{2}} \\ 
\mbox{}     &  &  &  &  & \vdots  &  &  &  &  & \\
\mbox{}G^{p^{e-1}} & = & H_{0}^{p^{e-1}} & \geq & H_{1}^{p^{e-1}} & \geq &
\cdots & \geq & H_{r}^{p^{e-1}} & = & G^{p^{e}}=\{1\}.
\end{array}$$
Omitting repetitions we get a chain where for each $k=0,\ldots ,e-1$ the length of the subchain between $G^{p^{k}}$ and $G^{p^{k+1}}$ has length equal to $\mbox{rank}(G^{p^{k}})$. Writing the groups in ascending order we get 
a chain of the form
    $$\{1\}=K_{0}<K_{1}<\cdots <K_{n}=G$$
where $|K_{i}|=p^{i}$ for $i=0,\ldots ,n$. We will now see how we can modify this
further. \\ \\
For $0\leq i\leq r-1$ we have $H_{i}=H_{i+1}\langle a_{i+1}\rangle$ and by
Lemma 2.1 we can choose our generators $a_{1},\ldots ,a_{r}$ such 
that $H_{i}^{p}=H_{i+1}^{p}$ if and only if $a_{i+1}^{p}=1$. Having done this we
can then move all the generators that are of order $p$ to the front of
the others (keeping the previous order unchanged otherwise) and we
still have that (1) gives us a powerfully central series. We can thus assume 
that for some $0\leq s\leq r$ we have $a_{1}^{p}=\cdots =a_{s}^{p}=1$ and that
$G^{p}=H_{s}^{p}>H_{s+1}^{p}>\cdots >H_{r}^{p}=G^{p^{2}}$. Notice that the rank of $G^{p}$ is the number of jumps and thus the number of generators that have order at least
$p^{2}$. We have $G^{p}=\{1\}$ if $s=r$ and otherwise $0\leq s<r$ and $\{1\}<G^{p}=H_{s}^{p}$. Suppose
we are in the latter situation. Using again Lemma 2.1 we see that for $s\leq j\leq r-1$ we have $H_{j}^{p^{2}}=H_{j+1}^{p^{2}}$ if and only if there exists
$x=y^{p}\in H_{j}^{p}\setminus H_{j+1}^{p}$ such that $x^{p}=1$. This happens if
and only if there is $y\in K_{j}\setminus K_{j+1}$ such that the order
of $y$ is $p^{2}$. We can thus choose our generators such that furthermore $H_{j}^{p^{2}}=K_{j+1}^{p^{2}}$ if and only if $a_{j+1}^{p^{2}}=1$. Notice again that the rank
of $G^{p^{2}}$ is the number of generators that have order at least $p^{3}$. Continuing in
this manner, considering next $H_{0}^{p^{3}}\geq H_{1}^{p^{3}}\geq \cdots 
\geq H_{r}^{p^{3}}=G^{p^{4}}$ and then $p^{4}$th powers and so on, we eventually arrive at a set of generators $a_{1},\ldots ,a_{r}$ with some specific properties. If for $0\leq i\leq r$
we let $s(i)$ be the number of generators of order $p^{i}$ then $|G^{p^{i-1}}/G^{p^{i}}|=p^{s(i)+s(i+1)+\cdots +s(e)}$. Using the fact that  $[G^{p^{i}},G]\leq G^{p^{i+1}}$ it is
easy to see inductively that every element in $G$ can be written of the form $a_{1}^{l_{1}}\cdots a_{r}^{l_{r}}$. We thus see from this that 
      $$G=\langle a_{1}\rangle\cdot \langle a_{2}\rangle\cdots
       \langle a_{r}\rangle$$
where 
 \begin{eqnarray*}
   |G| & = & |G/G^{p}|\cdot |G^{p}/G^{p^{2}}|\cdots |G^{p^{r-1}}/G^{p^{r}}| \\
       & = & p^{s(1)+\cdots + s(e)}p^{s(2)+\cdots + s(e)}\cdots p^{s(e)} \\
      & = & p^{s(1)}p^{2s(2)}\cdots p^{es(e)} \\
     & = & o(a_{1})\cdots o(a_{r})
\end{eqnarray*}
Notice that the number of generators of order $p^{i}$ is an invariant of the group, namely $s(i)=\mbox{rank\,}G^{p^{i}}-\mbox{rank\,}G^{p^{i+1}}$.
\begin{theo} Let $G$ be any powerfully nilpotent $p$-group of rank $r$,exponent $p^{e}$ and order $p^{n}$. Then we can choose our generators $a_{1},\ldots ,a_{r}$ such that $|G|=o(a_{1})\cdot
o(a_{2})\cdots o(a_{r})$ and such that for $H_{0}=G= \langle a_{1},\ldots ,a_{r}
\rangle G^{p}$, $H_{1}=\langle a_{2},\ldots ,a_{r}\rangle G^{p},\ldots ,H_{r}=G^{p}$
we get a powerfully central chain
$$\begin{array}{lllllllllll}
\mbox{}    G  & = & H_{0} & > & H_{1}  & > & \cdots & > & H_{r} & = & G^{p} \\
\mbox{}   G^{p} & = & H_{0}^{p} & \geq & H_{1}^{p} & \geq & \cdots & \geq & H_{r}^{p} & = &
   G^{p^{2}} \\ 
\mbox{}     &  &  &  &  & \vdots  &  &  &  &  & \\
\mbox{}G^{p^{e-1}} & = & H_{0}^{p^{e-1}} & \geq & H_{1}^{p^{e-1}} & \geq &
\cdots & \geq & H_{r}^{p^{e-1}} & = & G^{p^{e}}=\{1\}.
\end{array}$$
Moreover if we write the groups of this chain in ascending order without repetitions, then we get a powerfully central chain of the form
     $$\{1\}=K_{0}<K_{1}<\cdots <K_{n}=G$$
where $|K_{i}|=p^{i}$ for $i=1,\ldots, n$ and where the groups $G,G^{p},
\ldots ,G^{p^{e}}=\{1\}$ are included. Furthermore the generators satisfy relations of the form
$$
    [a_{i},a_{j}]=a_{1}^{m_{1}(i,j)}a_{2}^{m_{2}(i,j)}\cdots 
               a_{r}^{m_{r}(i,j)}$$
for $1\leq j<i\leq r$ where all the power indices are divisible by $p$ and where furthermore $p^{2}|m_{k}(i,j)$ when $k\leq i$. $G$ is then the largest finite $p$-group satisfying these relations and the relations
    $$a_{1}^{o(a_{1})}=\cdots =a_{r}^{o(a_{r})}.$$
The number of generators of any given order $p^{i}$ is also an invariant for the group $G$.
\end{theo}
{\bf Proof}\ \ Most of this has been proved already. The fact that $G$ satisfies
relations of the form given follows immediately from the fact that 
   $$G=H_{0}>H_{1}>\dots >H_{r}=G^{p}$$
is powerfully central. Now let $H=\langle b_{1},\ldots ,b_{r}\rangle$ be the largest finite $p$-group satisfying
the given relations in the variables $b_{1},\ldots ,b_{r}$. The relations then imply that $H$ is powerful and that the chain 
   $$\langle b_{1},\ldots ,b_{r}\rangle H^{p}>\langle b_{2},\ldots, b_{r}\rangle
H^{p}>\cdots >\langle b_{r}\rangle H^{p}>H^{p}$$
is powerfully central and thus $H$ is powerfully central by Proposition 1.1.
It follows that $H=\langle b_{1}\rangle \cdot \langle b_{2}\rangle \cdots
\langle b_{r}\rangle$ and thus $|H|\leq o(b_{1})\cdots o(b_{r})=o(a_{1})\cdots 
o(a_{r})=|G|$. As $G$ is clearly a homomorphic image of $H$ it follows that
$|H|=|G|$ and thus $H\cong G$. $\Box$. \\ \\
{\bf Definition}. Consider a presentation of a group $G=\langle b_{1},\ldots ,
b_{r}\rangle$ satisfying relations of the form 
$$
    [b_{i},b_{j}]=b_{1}^{m_{1}(i,j)}b_{2}^{m_{2}(i,j)}\cdots 
               b_{r}^{m_{r}(i,j)}$$
for $1\leq j<i\leq r$ where all the power indices are divisible by $p$ and where furthermore $p^{2}|m_{k}(i,j)$ when $k\leq i$. To these we add relations of the form 
    $$b_{1}^{p^{e_{1}}}=1,\ b_{2}^{p^{e_{2}}}=1 \cdots ,b_{r}^{p^{e_{r}}}=1$$
for some positive integers $e_{1},\ldots ,e_{r}$. We call such a presentation 
a {\it powerfully nilpotent presentation}. \\ \\
{\bf Remark}. We know from Theorem 2.7 that every powerfully nilpotent $p$-group has a presentation of this form. \\ \\
{\bf Definition}. We say that a powerfully nilpotent presentation is consistent if the largest finite $p$-group $H$ satisfying the relations has order
$p^{e_{1}}p^{e_{2}}\cdots p^{e_{r}}$. 
\section{Bounding the exponent in terms of the powerful coclass} 
\begin{lemm} Let $G$ be any powerful $p$-group and let $a,b\in G$. For any integer $k\geq 0$ we have 
     $$[a^{p^{k}},b]=[a,b]^{p^{k}}=[a,b^{p^{k}}]$$
modulo $[G^{p^{k+1}},G]$. 
\end{lemm} 
{\bf Proof}. \  \  First suppose that $p$ is odd. We prove this by induction on $k$. For $k=0$ this is clear. Now suppose $k\geq 1$ and that the result holds for smaller values of $k$. Using the Interchange Lemma we observe
that $[G^{p^{k-1}},G,G]^{p}=[[G^{p^{k-1}},G]^{p},G]$, $[G^{p^{k-1}},G,G,G]$, $[G^{p^{k}},G,G]$ and $[G^{p^{k}},G]^{p}$
are all contained in $[G^{p^{k+1}},G]$. We will be using this fact in the following calculations. By the induction hypothesis we have
that $[a^{p^{k-1}},b]=[a,b]^{p^{k-1}}u$ for some $u\in [G^{p^{k}},G]$. Let $x=a^{p^{k-1}}$. Then, calculating modulo $[G^{p^{k+1}},G]$, we have
\begin{eqnarray*}
     [a^{p^{k}},b] & = & [x^{p},b] \\
                 & = & [x,b]^{p}[x,b,x]^{p\choose 2} \\
                 & = & [a^{p^{k-1}},b]^{p} \\
                 & = & ([a,b]^{p^{k-1}}u)^{p} \\
                 & = & [a,b]^{p^{k}}u^{p} \\
                 & = & [a,b]^{p^{k}}.
\end{eqnarray*}
This proves the induction step and thus completes the proof when $p$ is odd. When $p=2$, the induction step is easier as $[G,G]\leq G^{4}$ and thus $[G^{p^{k-1}},G,G]\leq [G^{p^{k+1}},G]$. 
Using this fact, the calculations as above show that, modulo $[G^{p^{k+1}},G]$, we have $[a^{p^{k}},b]=[a,b]^{p^{k}}$. $\Box$  \\ \\
{\bf Remark}. As $[G^{p^{k+1}},G]\leq G^{p^{k+2}}$ we get in particular that $[a^{p^{k}},b]=[a,b]^{p^{k}}=[a,b^{p^{k}}]$ modulo $G^{p^{k+2}}$. 
\begin{lemm} Let $G$ be a powerful $p$-group and suppose $G=\langle a_{1},\ldots ,a_{r}\rangle$ and that 
$G^{p^{k}}=\langle a_{i_{1}}^{p^{k}},\ldots ,a_{i_{s}}^{p^{k}}\rangle$ for some
$1\leq i_{1}<i_{2}<\cdots < i_{s}\leq r$. Then 
     $$[\langle a_{i_{j}}^{p^{k}}\rangle ,G][G^{p^{k+1}},G]=
     \langle [a_{i_{j}}^{p^{k}},a_{i_{1}}],\ldots ,[a_{i_{j}}^{p^{k}},
a_{i_{s}}]\rangle [G^{p^{k+1}},G].$$
\end{lemm}
{\bf Proof}\ \ It clearly suffices to show that for $m\in \{1,\ldots ,r\}$ we
have
      $$[a_{i_{j}}^{p^{k}},a_{m}]\in \langle [a_{i_{j}}^{p^{k}},a_{i_{1}}],
\ldots ,[a_{i_{j}}^{p^{k}},a_{i_{s}}]\rangle [G^{p^{k+1}},G].$$
Now $a_{m}^{p^{k}}=a_{i_{1}}^{p^{k}e_{1}}a_{i_{2}}^{p^{k}e_{2}}\cdots a_{i_{s}}^{p^{k}e_{s}}$ for some integers $e_{1},\ldots ,e_{s}$. Thus, calculating modulo $[G^{p^{k+1}},G]$ and using Lemma 3.1, we have 
\begin{eqnarray*}
    [a_{i_{j}}^{p^{k}},a_{m}] & = & [a_{i_{j}},a_{m}^{p^{k}}] \\
      & = & [a_{i_{j}},a_{i_{1}}^{p^{k}e_{1}}\cdots a_{i_{s}}^{p^{k}e_{s}}] \\
     & = & [a_{i_{j}},a_{i_{1}}^{p^{k}}]^{e_{1}}\cdots [a_{i_{j}},a_{i_{s}}^{p^{k}}]^{e_{s}} \\
     & = & [a_{i_{j}}^{p^{k}},a_{i_{1}}]^{e_{1}}\cdots [a_{i_{j}}^{p^{k}},a_{i_{s}}]^{e_{s}}.
\end{eqnarray*}
This finishes the proof. $\Box$ \\
\begin{coro} Let $G$ be a powerful $p$-group. If the rank of $G^{p^{k}}$ is $1$ then $G^{p^{k}}\leq Z(G)$. 
\end{coro} 
{\bf Proof}. Suppose $G^{p^{k}}=\langle c^{p^{k}}\rangle$. By Lemma 3.2 we then have 
      $$[G^{p^{k}},G]=\langle [c^{p^{k}},c]\rangle [G^{p^{k+1}},G]=
			         [G^{p^{k+1}},G].$$
But then $[G^{p^{k}},G]=[G^{p^{k}},G]^{p}=[G^{p^{k}},G]^{p^{2}}
=\cdots =\{1\}$. $\Box$ \\
\begin{prop} Let $G$ be a powerfully nilpotent group and suppose that $G^{p^{k}}$ has rank $s\geq 2$. There exists a chain
    $$G^{p^{k}}=H_{0}>H_{1}>\ldots >H_{s-1}=G^{p^{k+1}}$$
that is powerfully centralized by $G$.
\end{prop}
{\bf Proof}\ \ Suppose $G$ has rank $r$. By Theorem 2.7, we can pick our generator $a_{1},\ldots ,a_{r}$ such that for
$K_{0}=\langle a_{1},\ldots ,a_{r}\rangle G^{p}, K_{1}=\langle a_{2},\ldots ,a_{r}\rangle G^{p},\ldots ,K_{r}=G^{p}$ we have 
a chain that is powerfully centralized by $G$. Furthermore the chain \\ 
$$G^{p^{k}}=K_{0}^{p^{k}}\geq K_{1}^{p^{k}}\geq \cdots \geq K_{r}^{p^{k}}=G^{p^{k+1}}$$
is also powerfully centralized by $G$.
Omitting repetitions we
get a chain of the form
 $$G^{p^{k}}=\langle a_{i_{1}}^{p^{k}},\ldots ,
a_{i_{s}}^{p^{k}}\rangle G^{p^{k+1}} > \langle a_{i_{2}}^{p^{k}},\ldots,
a_{i_{s}}^{p^{k}}\rangle G^{p^{k+1}}>\ldots > \langle a_{i_{s}}^{p^{k}}\rangle G^{p^{k+1}}>G^{p^{k+1}}$$
for some $1\leq i_{1}<i_{2}<\cdots <i_{s}\leq r$. We can do better than this.
Using Lemmas 3.1 and 3.2, we have 
\begin{eqnarray*}
    [\langle a_{i_{1}}^{p^{k}}\rangle, G]G^{p^{k+2}} & = & 
      \langle [a_{i_{1}}^{p^{k}},a_{i_{2}}],\ldots ,
         [a_{i_{1}}^{p^{k}},a_{i_{s}}]\rangle G^{p^{k+2}} \\
     & = &  \langle [a_{i_{1}},a_{i_{2}}^{p^{k}}],\ldots [a_{i_{1}},a_{i_{s}}^{p^{k}}]\rangle G^{p^{k+2}}.
\end{eqnarray*}
As $[G^{p^{k+1}},G]\leq G^{p^{k+2}}$, this implies that 
$$[G^{p^{k}},G]\leq [\langle a_{i_{2}}^{p^{k}},\ldots ,a_{i_{s}}^{p^{k}}\rangle G^{p^{k+1}},G]G^{p^{k+2}}\leq
(\langle a_{i_{3}}^{p^{k}},\ldots ,a_{i_{s}}^{p^{k}}\rangle G^{p^{k+1}})^{p}.$$ 
We thus have that 
 $$G^{p^{k}}=\langle a_{i_{1}}^{p^{k}},\ldots ,
a_{i_{s}}^{p^{k}}\rangle G^{p^{k+1}} > \langle a_{i_{3}}^{p^{k}},\ldots,
a_{i_{s}}^{p^{k}}\rangle G^{p^{k+1}}>
\langle a_{i_{4}}^{p^{k}},\ldots ,a_{i_{s}}^{p^{k}}\rangle G^{p^{k+1}}>\ldots > \langle a_{i_{s}}^{p^{k}}\rangle G^{p^{k+1}}>G^{p^{k+1}}$$     
is powerfully central. $\Box$ \\ \\
\begin{theo} Let $G$ be a powerfully nilpotent group of order $p^{n}$, powerful class $c$ and exponent $p^{e}$. Then 
             $$e\leq n-c+1.$$
\end{theo}
{\bf Proof}\ \ This is easy to see when $G$ is cyclic so we can assume that the rank of $G$ is
at least $2$. Let $k$ be the largest  non-negative  integer such that the rank of
$G^{p^{k}}$ is greater than or equal to  $2$. Let $r_{i}$ be the rank of $G^{p^{i}}$ for 
$i=0,1,\ldots ,k$ and let $p^{n_{0}}=|G^{p^{k+1}}|$. Notice then that
\begin{eqnarray*}
    e & = & k+1+n_{0} \\
   n & = & r_{0}+r_{1}+\cdots +r_{k}+n_{0}.
\end{eqnarray*}
By Proposition 3.4 there exists for each $0\leq j\leq k$ a descending chain 
    $$G^{p^{j}}=H_{0}>H_{1}>\ldots >H_{r_{j}-1}=G^{p^{j+1}}$$
that is powerfully centralized by $G$. 
Adding up for $j=0,1,\ldots ,k$ and using the fact from Corollary  3.3  that $G^{p^{k+1}}\leq Z(G)$, we get a central chain of
total length $(r_{0}-1)+(r_{1}-1)+\cdots +(r_{k}-1)+1$. Hence 
   $$c\leq (r_{0}-1)+(r_{1}-1)+\cdots + (r_{k}-1)+1.$$
We conclude that 
\begin{eqnarray*}
   n-c & \geq  & n-[(r_{0}-1)+\cdots +(r_{k}-1)+1] \\
      & = & (r_{0}+r_{1}+\cdots +r_{k}+n_{0})-[(r_{0}-1)+\cdots +(r_{k}-1)+1] \\
     & = & k+n_{0} \\
     & = & e-1.
\end{eqnarray*}
Hence $e\leq n-c+1$. $\Box$ \\ \\
As a corollary we get one of the central results of this paper. \\
\begin{theo} For each prime $p$ and non-negative integer $d$, there are finitely many powerfully nilpotent $p$-groups of powerful coclass $d$.
\end{theo}
{\bf Proof}\ \ Let $G$ be a powerfully nilpotent group of order $p^{n}$, rank $r$ and exponent $e$. Then $n\leq re\leq (d+1)(d+1)$ and thus the order of $G$ is bounded by the powerful coclass. $\Box$
\section{Groups with maximal tail}
{\bf Definition}. Let $G$ be a powerfully nilpotent $p$-group and let $k$ be
the largest non-negative integer such that 
   $$p=|Z(G)^{p}|=|\frac{\hat{Z}_{2}(G)^{p}}{\hat{Z}_{1}(G)^{p}}|=\cdots 
      =|\frac{\hat{Z}_{k}(G)^{p}}{\hat{Z}_{k-1}(G)^{p}}|.$$
we refer to $\hat{Z}_{k}(G)^{p}$ as the {\it tail} of $G$ and $k$ as the {\it length of the tail}. \\ \\
{\bf Remark}. If $G$ has a tail of length $k$ then $G,
G/\hat{Z}_{1}(G)^{p}, G/\hat{Z}_{2}(G)^{p},\ldots ,G/\hat{Z}_{k}(G)^{p}$ all have the same powerful coclass. \\ \\
Now let $G$ be any powerfully nilpotent $p$-group of rank $r$ and exponent $p^{e}$. By Theorem 2.7 we can find $a_{1},a_{2},\ldots ,a_{r}\in G$
such that for 
$$K_{0}=\langle a_{1},a_{2}\ldots, a_{r}\rangle G^{p},
K_{1}=\langle a_{2},\ldots ,a_{r}\rangle G^{p},\cdots ,K_{r}=G^{p}$$
we have that the chain 
$$\begin{array}{lllllllllll}
\mbox{}    G  & = & K_{0} & \geq  & K_{1}  & \geq  & \cdots & \geq  & K_{r} & = & G^{p} \\
\mbox{}   G^{p} & = & K_{0}^{p} & \geq & K_{1}^{p} & \geq & \cdots & \geq & K_{r}^{p} & = &
   G^{p^{2}} \\ 
\mbox{}     &  &  &  &  & \vdots  &  &  &  &  & \\
\mbox{}G^{p^{e-1}} & = & K_{0}^{p^{e-1}} & \geq & K_{1}^{p^{e-1}} & \geq &
\cdots & \geq & K_{r}^{p^{e-1}} & = & G^{p^{e}}=\{1\},
\end{array}$$
is powerfully central in $G$.  
\begin{lemm}
Rewrite the chain in ascending order without repetitions.
Then suppose the chain up to and including $G^{p}$ is
    $$\{1\}=M_{0}<M_{1}<\ldots <M_{t}=G^{p}.$$
We have that $M_{j}\leq \hat{Z}_{j}(G)^{p}$ for $j=0,\ldots ,t$. Also if the tail of $G$ is $\hat{Z}_{k}(G)^{p}$ then 
$M_{i}=\hat{Z}_{i}(G)^{p}$ for $i=0,\ldots ,k$. 
\end{lemm}
{\bf Proof}\ \ We prove the first part by induction on $0\leq j\leq t$. This is obvious when $j=0$. Now suppose that $j\geq 1$ and that the result holds for
smaller values of $j$. Let $q$ be the largest and then, for that $q$, $i$ 
be the largest such that $M_{j}=K_{i}^{p^{q}}$. Then $0\leq i\leq r-1$ and
$K_{i+1}^{p^{q}}=M_{j-1}$. Thus 
   $$[K_{i}^{p^{q-1}},G]\leq (K_{i+1}^{p^{q-1}})^{p}=M_{j-1}\leq \hat{Z}_{j-1}(G)^{p}$$
by the induction hypothesis. Hence $K_{i}^{p^{q-1}}\leq \hat{Z}_{j}(G)$ and thus
$M_{j}=K_{i}^{p^{q}}\leq \hat{Z}_{j}(G)^{p}$. This finishes the inductive proof. The 2nd part follows from the 1st part and the fact that $|\hat{Z}_{i}(G)^{p}|=p^{i}$ for $i=0,\ldots ,k$. $\Box$
\mbox{}\\ 
\begin{prop} Let $G$ be a powerfully nilpotent $p$-group with a tail of length $k$. Suppose $G^{p^{i+1}}\leq \hat{Z}_{k}(G)^{p}$ for some non-negative integer $i$. If $\mbox{rank\,}(G^{p^{i}})\geq 2$ we have
    $$\mbox{rank\,}(G^{p^{i}})>\mbox{rank\,}(G^{p^{i+1}}).$$
\end{prop}
{\bf Proof}\ \ Suppose the rank of $G^{p^{i}}$ is $s$. We know that the following chain is powerfully centralized by $G$
 $$G^{p^{i}}=K_{0}^{p^{i}}\geq K_{1}^{p^{i}}\geq \cdots
      \geq K_{r}^{p^{i}}=G^{p^{i+1}}.$$
Omitting repetitions we get a chain
  $$G^{p^{i}}=\langle a_{i_{1}}^{p^{i}},\ldots ,a_{i_{s}}^{p^{i}}
	\rangle G^{p^{i+1}}>\langle a_{i_{2}}^{p^{i}},\ldots ,
	   a_{i_{s}}^{p^{i}
		}\rangle G^{p^{i+1}}> \cdots >\langle a_{i_{s}}^{p^{i}}\rangle G^{p^{i+1}}>G^{p^{i+1}}$$
with $1\leq i_{1}<i_{2}<\cdots <i_{s}\leq r$.
As we saw in the proof of Proposition 3.4 we can omit the 2nd term and we still have a chain powerfully centralized by $G$. All the terms are then in particular powerfully embedded in $G$. If $G^{p^{i}}\geq E>F\geq G^{p^{i+1}}$ where $E,F$ are consecutive terms in the chain, then $[E^{p},G]=[E,G]^{p}\leq (F^{p})^{p}$. Thus we get the chain 
$$\begin{array}{lllllllll}

G^{p^{i}}=\langle a_{i_{1}}^{p^{i}},\ldots ,a_{i_{s}}^{p^{i}}
	\rangle G^{p^{i+1}} & > & \langle a_{i_{3}}^{p^{i}},\ldots ,
	   a_{i_{s}}^{p^{i}
		}\rangle G^{p^{i+1}} & > & \cdots >\langle a_{i_{s}}^{p^{i}}\rangle G^{p^{i+1}} & > & G^{p^{i+1}} \\
G^{p^{i+1}}=\langle a_{i_{1}}^{p^{i+1}},\ldots ,a_{i_{s}}^{p^{i+1}}
	\rangle G^{p^{i+2}} & \geq & \langle a_{i_{3}}^{p^{i+1}},\ldots ,
	   a_{i_{s}}^{p^{i+1}
		}\rangle G^{p^{i+2}} & \geq & \cdots \geq \langle a_{i_{s}}^{p^{i+1}}\rangle G^{p^{i+2}} & \geq & G^{p^{i+2}}
\end{array}$$
that is powerfully centralized by $G$. Let $E=G^{p^{i}}$ and $F=\langle a_{i_{3}}^{p^{i}},\ldots
,a_{i_{s}}^{p^{i}}\rangle G^{p^{i+1}}$. By Lemma 4.1
we know that $F^{p}=K_{i_{3}-1}^{p^{i+1}}$ is equal to some $\hat{Z}_{t}(G)^{p}$ where $0\leq t\leq k$. 
As $[E,G]\leq F^{p}=\hat{Z}_{t}(G)^{p}$, it follows
that $E\leq \hat{Z}_{t+1}(G)$ and thus $E^{p}\leq
\hat{Z}_{t+1}(G)^{p}$. As $E^{p}=G^{p^{i+1}}\leq 
\hat{Z}_{k}(G)^{p}$ we know from Lemma 4.1 that 
$[E^{p}:F^{p}]\leq [\hat{Z}_{t+1}(G)^{p}:\hat{Z}_{t}(G)^{p}]=p$ and it follows that the rank of $G^{p^{i+1}}$ is at most $s-1$. $\Box$ \\ 
\begin{theo} Let $G$ be a powerfully nilpotent group of rank $r\geq 2$ and let $f$ be the largest non-negative integer such that $G^{p^{f}}$ has
rank at least $2$. Let $i$ be the smallest non-negative integer such that $G^{p^{i+1}}$ is contained in the tail of $G$. Then \\ \\
(a) $\mbox{rank}(G^{p^{i}})>\mbox{rank}(G^{p^{i+1}})>\cdots >\mbox{rank}(G^{p^{f+1}})$. \\ 
(b) The length of the tail of $G$ is at most $1+\frac{r(r-1)}{2}$. 
\end{theo}
{\bf Proof}\ \ (a) follows directly from Proposition 4.2. For (b) notice that the scenario where the length of the tail would potentially be largest is 
when the tail
is equal to $G^{p}$ and where $\mbox{rank}(G)=r,
\mbox{rank}(G^{p})=r-1, \ldots, \mbox{rank}(G^{p^{r-1}})=1$. By Corollary 3.3 we then have that $G^{p^{r-1}}\leq Z(G)$. For there to be a tail of length greater than $0$ we then need $|Z(G)^{p}|=p$ and thus $|G^{p^{r-1}}|\leq p^{2}$. Thus the size of the tail can't be longer than
    $$|G^{p}|=|G^{p}/G^{p^{2}}|\cdot|G^{p^{2}}/G^{p^{3}}|\cdots |G^{p^{r-2}}/G^{p^{r-1}}|\cdot |G^{p^{r-1}}|\leq p^{(r-1)+(r-2)+\cdots +2+2}.$$
This finishes the proof. $\Box$ \\ \\
{\bf Definition}. Let $G$ be a powerfully nilpotent $p$-group. We say that $G$ has {\it maximal tail} if the tail of $G$ is $G^{p}$. \\ \\
{\bf Definition}. Let $G$ be a powerfully nilpotent $p$-group. We say that a subgroup $H$ is 
powerfully hypercentral if there exists a descending chain of subgroups
     $$H=H_{0}>H_{1}>\ldots >H_{n}=\{1\}$$
such that $[H_{i},G]\leq H_{i+1}^{p}$ for $i=0,\ldots ,n-1$. \\
\begin{theo} Let $G$ be a powerfully nilpotent group
of rank $r\geq 2$ that has a maximal tail. Suppose $G$ has order $p^{n}$, powerful class $c$ and exponent
$p^{e}$. Let $t$ be the length of the tail. \\ \\
(a) We have that $t=n-r$ and $c-1\leq t\leq c$. It follows also that $n-c\leq r\leq n-c+1$. 
(b) We have $t\leq 1+\frac{r(r-1)}{2}$ and
$n\leq 1+\frac{r(r+1)}{2}$. \\ 
(c) We have $\mbox{rank}(G)>\mbox{rank}(G^{p})>\cdots 
>\mbox{rank}(G^{p^{e-2}})$. \\
(d) We have $G/\hat{Z}_{i}(G)^{p}$ has maximal tail for all non-negative integers $i$. \\ 
(e) If $H\leq G$ is powerfully hypercentral, then $H^{p}=\hat{Z}_{i}(G)^{p}$ for some integer $i\geq 0$. 
\end{theo}
{\bf Proof} (a) As $|G^{p}|=p^{n-r}$, it is clear that $t=n-r$.  From Lemma 2.5 we know that $\hat{Z}_{c-2}(G)^{p}<\hat{Z}_{c-1}(G)^{p}\leq \hat{Z}_{c}(G)^{p}=G^{p}$. As the tail is $G^{p}$ we see from this that $c-1\leq t\leq c$.  \\ \\
(b) The first inequality follows from Theorem 4.3 and the 2nd one follows from $|G|=|G^{p}|\cdot |G/G^{p}|=p^{t}p^{r}$. \\ \\
(c) If $\mbox{rank}(G^{p^{e-1}})\geq 2$, this follows from Theorem 4.3. Now suppose that $\mbox{rank}(G^{p^{e-1}})=1$. Let $f$ be the largest non-negative integer such that $G^{p^{f}}$ has rank at least $2$. Then $G^{p^{f+1}}$ is cyclic and non-trivial and by Corollary 3.3 we then know that $G^{p^{f+1}}\leq
Z(G)$. As $|Z(G)^{p}|=p$ it follows that $G^{p^{f+3}}=\{1\}$ and thus $e\leq f+3$ that implies that
$e-2\leq f+1$. The result now follows from Theorem 4.3 and the fact that $e-2$ is the smallest integer
such that $G^{p^{e-2}}$ has rank $1$. \\ \\
(d) This follows straight from the fact that $\hat{Z}_{i}(G)^{p}\leq G^{p}$ and that $G$ has maximal tail. \\ \\
(e) If $H^{p}=\{1\}$ then $H^{p}=\hat{Z}_{0}(G)^{p}$. Thus suppose that $H^{p}\not =\{1\}$. By our assumptions there is a descending chain
    $$H=H_{0}>H_{1}>\ldots >H_{n}=\{1\}$$
where $[H_{i},G]\leq H_{i+1}^{p}$ for $i=0,\ldots,
n-1$. Let $j$ be smallest such that $H_{j}^{p}=\{1\}$.  Notice that $j\geq 1$. As $[H_{j-1},G]\leq
H_{j}^{p}=\{1\}$ we have $H_{j-1}\leq Z(G)$. Then
$\{1\}<H_{j-1}^{p}\leq Z(G)^{p}$ and as $|Z(G)^{p}|=p$, it follows that $H_{j-1}^{p}=Z(G)^{p}$. Replacing $H,G$ by $H/Z(G)^{p}, G/Z(G)^{p}$, we see inductively that $H^{p}=\hat{Z}_{j}(G)^{p}$ for some
$j$. $\Box$ \\ \\
{\bf Remark}. (1) We know from Theorem 4.4 that if $G$ has maximal tail then $n-c\leq r\leq n-c+1$.  Let $G$ be a powerfully nilpotent group of rank $r\geq 2$, powerful class $c$, order $p^{n}$ and exponent $p^{e}$ where $r=n-c+1$. Notice then $|G^{p}|=
p^{n-r}=p^{c-1}$. From Lemma 2.5 we know that
   $$\{1\}=\hat{Z}_{0}(G)^{p}<\hat{Z}_{1}(G)^{p}<
	\cdots <\hat{Z}_{c-1}(G)^{p}.$$
Hence $Z_{c-1}(G)^{p}=G^{p}$ and $G$ has maximal tail of length $c-1$. \\ \\
If on the other hand $r=n-c$ and $G$ has maximal tail, then the length of the tail must be $c$ and this happens if and only if $\hat{Z}_{c-1}(G)^{p}<
\hat{Z}_{c}(G)^{p}$. \\ \\
(2) We know from Theorem 4.4 that if $G$ has maximal tail, then the length of the tail is $t\leq 1+r(r-1)/2$. We also know from the proof of Theorem 
4.3 that in order for the upper bound to be attained, we need $\mbox{rank}(G)=r, \mbox{rank}(G^{p})=r-1,\ldots ,\mbox{rank}(G^{p^{r-2}})=2$ and
that $G^{p^{r-1}}$ is cyclic of order $p^{2}$. In particular we must have that $e=r+1$. As $n-c\leq r\leq n-c+1$ and $e\leq n-c+1$ this can only
happen if 
\begin{equation}
                                      r=n-c\mbox{\ and\ }e=n-c+1.
\end{equation}
Conversely,  $G$ has maximal tail of length $1+r(r-1)/2$   if $\hat{Z}_{c-1}(G)^{p}<\hat{Z}_{c}(G)^{p}$ and (2) holds. To see this notice first that $r=n-c$ and $\hat{Z}_{c-1}(G)^{p}<\hat{Z}_{c}(G)^{p}$ implies that $|\hat{Z}_{c}(G)^{p}|\geq p^{c}=p^{n-r}=|G^{p}|$ and thus $G^{p}=\hat{Z}_{c}(G)^{p}$. Hence
this implies that $G$ has maximal tail.  From Theorem 4.4 we then have that 
   $$\mbox{rank}(G)=r,\mbox{rank}(G^{p})=r-1,\ldots ,\mbox{rank}(G^{p^{r-1}})=\mbox{rank}(G^{p^{e-2}})=1$$
 and thus, as $e=r+1$, $G^{p^{r-1}}$ is cyclic of order $p^{2}$. It follows then that $|G^{p}|=p^{1+r(r-1)/2}$ and thus $t=1+r(r-1)/2$.   We conjecture that for each $r\geq 1$ there exist a group $G$ with maximal tail of length $t=1+r(r-1)/2$. The following examples give such examples for $r\leq 5$.  \\ \\
{\bf Example 1} ($r\leq 4$). Let $G(r)=\langle x,a_{1},\ldots ,a_{r-1}\rangle$ where the following relations hold
 $$\begin{array}{l}
   x^{p^{r+1}}=1,\, a_{1}^{p}=a_{2}^{p^{2}}=\cdots =a_{r-1}^{p^{r-1}}=1 ,\\
\mbox{}[a_{1},x]=a_{2}^{p},\ [a_{2},x]=a_{3}^{p},\ldots ,[a_{r-2},x]=a_{r-1}^{p},\, [a_{r-1},x]=x^{p^{2}}.
\end{array}$$
\mbox{}\\ \\
Unfortunately, this recipe does not work for higher $r$ without modification, although we believe that a modification in line with the one given in the following example will work for all $r$. \\ \\
{\bf Example 2}. Let $G=\langle x,a,b,c,d\rangle$ with relations
 $$\begin{array}{l}
    x^{15625}=1,\, a^{5}=b^{25}=c^{125}=d^{625}=1, \\
\mbox{}[a,x]=b^{5},\, [b,x]=c^{5},\, [c,x]=d^{5},\, [d,x]=x^{25},\, [c,d]=c^{25}d^{375}.
\end{array}$$
%
%Suppose furthermore that $e=n-c+1=r$. This can
%only happen if 
%
 %  $$\mbox{rank}(G)=r, \mbox{rank}(G^{p})=r-1,
%	\cdots ,\mbox{rank}(G^{p^{r-1}})=1.$$
%
%Thus we must have that the tail is $c-1$
%where 
%
 %   $$p^{c-1}=|G^{p}|=p^{(r-1)+(r-2)+\cdots +1}.$$
%
%Thus we get
%
%\begin{eqnarray*}
%
%   c & = & 1+\frac{r(r-1)}{2} \\
%
  % n & = & c-1+r=\frac{r(r+1)}{2} \\
%
 %  e & = & r \\
%
 % n-c & = & r-1.
%
%\end{eqnarray*}
%
\section{Growth of PN groups of exponent $p^{2}$}
Let $p$ be a given odd prime. In this section we describe the growth of powerfully nilpotent $p$-groups of exponent $p^{2}$  in terms of the order. \\  \\
Let $y,x$ be fixed non-negative integers. Let $G=\langle a_{1}, a_{2}, \ldots ,a_{y+x}\rangle$ be a powerfully nilpotent group of rank $r=y+x$ and order
$p^{n}$ where $n=y+2x$. Here we are assuming furthermore that $o(a_{1})=\cdots =o(a_{y})=p$, $o(a_{y+1})=\cdots =o(a_{y+x})=p^{2}$ and that 
these satisfy a powerfully nilpotent presentation 
    $$[a_{i},a_{j}]=a_{i+1}^{p\alpha_{i+1}(i,j)}a_{i+2}^{p\alpha_{i+2}(i,j)}\cdots a_{y+x}^{p\alpha_{y+x}(i,j)}$$
for $1\leq j<i\leq y+x$ with $0\leq \alpha_{k}(i,j)\leq p-1$. \\ \\
Taking into account that $a_{1}^{p}=\ldots =a_{y}^{p}=1$ and that $y=n-2x$, we see that for $1\leq j<i\leq y=n-2x$ there are $p^{x}$ choices for $[a_{i},a_{j}]$. Notice that there are ${y\choose 2}={n-2x\choose 2}$ such pairs $(i,j)$.  For $y+1\leq i\leq y+x$ and $1\leq j<i$ there are $p^{x+y-i}$ possible values for $[a_{i},a_{j}]$ and that for each such $i$ there are $i-1$ such pairs $(i,j)$. Thus for a fixed $n$ the total number of powerfully nilpotent presentations is $p^{h(x)}$ where 
\begin{eqnarray*}
     h(x) & = & {y\choose 2}x+y(x-1)+(y+1)(x-2)+\cdots + (y+x-2)\cdot 1 \\
         & = & \frac{x}{6}(7x^{2}-9(n-1)x+3n^{2}-6n+2)
\end{eqnarray*}
and thus 
\begin{eqnarray*}
  h'(x)      & = & \frac{21}{6}[(x-\frac{3}{7}(n-1))^{2}-(\frac{2}{49}(n-1)^{2}+1/21)].
\end{eqnarray*}
The two roots of $h'(x)$ are $\frac{3}{7}(n-1)\pm  \sqrt{\frac{2}{49}(n-1)^{2}+1/21}$. Notice that $0\leq x\leq n/2$ and it is not difficult to see that for sufficiently large $n$
        $$0< \frac{3}{7}(n-1)-\sqrt{\frac{2}{49}(n-1)^{2}+1/21} < n/2 < \frac{3}{7}(n-1)+ \sqrt{\frac{2}{49}(n-1)^{2}+1/21}.$$
It follows from this that for sufficiently large $n$, the cubic $h(x)$ takes it's maximal value when $x=x_{max}=\frac{3}{7}(n-1)-\sqrt{\frac{2}{49}(n-1)^{2}+1/21}$ and where
the integer $x(n)$, that gives the largest value for $h(x)$, satisfies $x_{max}-1\leq x(n)\leq x_{max}+1$. Thus we have that
$x(n)/n$ tends to $(3-\sqrt{2})/7$ as $n$ tends to infinity. Thus $h(x(n))/n^{3}=\frac{x(n)}{6n}(7x(n)^{2}/n^{2}-9(1-1/n)x(n)/n+3-6/n+2/n^{2})$ tends to 
      $$\frac{3-\sqrt{2}}{6\cdot 7}(7(\frac{3-\sqrt{2}}{7})^{2}-9\frac{3-\sqrt{2}}{7}+3)=\frac{9+4\sqrt{2}}{294}$$
as $n$ tends to infinity. \\ \\
Let $n$ be fixed. For any $0\leq x\leq \lfloor\frac{n}{2}\rfloor$, let ${\mathcal P}(n,x)$ be the collection of all powerfully nilpotent presentations as above. It is not difficult to see that those presentations are consistent and thus the resulting group is of exponent $p^{2}$, order $p^{n}$ and rank $n-x$. Furthermore $a_{1}^{p}=\cdots =a_{n-2x}^{p}=1$ and $a_{n-2x+1}^{p^{2}}=\cdots =a_{n-x}^{p^{2}}=1$. We have just seen that if $x(n)$ is chosen
such that number of presentations is maximal then 
                         $$|{\mathcal P}(n,x(n))|=p^{\alpha n^{3}+o(n^{3})}$$
where $\alpha=\frac{9+4\sqrt{2}}{294}$. Let ${\mathcal P}_{n}$ be the total number of the powerfully nilpotent presentations where $0\leq x\leq \lfloor\frac{n}{2}\rfloor$. Then 
                  $${\mathcal P}_{n}={\mathcal P}(n,0)\cup {\mathcal P}(n,1)\cup \cdots \cup {\mathcal P}(n,\lfloor n/2\rfloor)$$
and thus
          $$p^{\alpha n^{3}+o(n^{3})}=|{\mathcal P}(n,x(n))|\leq |{\mathcal P}_{n}|=|{\mathcal P}(n,0)|+\cdots +|{\mathcal P}(n,\lfloor n/2\rfloor)|\leq n\cdot |{\mathcal P}(n,x(n))|=p^{\alpha n^{3}+o(n^{3})}.$$
This shows that 
                            $$|{\mathcal P}_{n}|=p^{\alpha n^{3}+o(n^{3})}.$$
We want to show that this also gives us the growth of powerfully nilpotent groups of exponent $p^{2}$ with respect to the order $p^{n}$. Clearly $p^{\alpha n^{3}+o(n^{2})}$ gives
us an upper bound. We want to show that this is also a lower bound. Let $x=x(n)$ be as above. Let $a_{1},\ldots ,a_{n-x}$ be a set of generators for a powerfully nilpotent group $G$  where
$a_{1}^{p}=\cdots =a_{n-2x}^{p}=1$ and $a_{n-2x+1}^{p^{2}}=\cdots =a_{n-x}^{p^{2}}=1$. Notice that $\langle a_{1},\ldots ,a_{n-2x}\rangle G^{p}=G[p]=\{g\in G:\,g^{p}=1\}$ and this is thus
a characteristic subgroup of $G$. It will be useful to consider a larger class of presentations for powerfully nilpotent groups of exponent $p^{2}$ and order $p^{n}$ where we still require
$o(a_{1})=\cdots =o(a_{n-2x})=p$ and $o(a_{n-2x+1})=\cdots o(a_{n-x})=p^{2}$. We let $Q(n,x)=Q(n,x(n))$ be the collection of all presentations with additional commutator relations 
                $$[a_{i},a_{j}]=a_{1}^{p\alpha_{1}(i,j)}\cdots a_{n-x}^{p\alpha_{n-x}(i,j)}.$$
The presentation will be included in $Q(n,x)$ provided the resulting group is a powerfully nilpotent group $G$ of exponent $p^{2}$ and order $p^{n}$.  
Notice then that $G^{p}\leq Z(G)$ and as a result the commutator relations above only depend on the cosets $\bar{a}_{1}=a_{1}G^{p},\ldots ,\bar{a}_{n-x}=a_{n-x}G^{p}$ and not on the exact values of $a_{1},\ldots ,a_{n-x}$. Consider the vector space $V=G/G^{p}={\mathbb Z}_{p}\bar{a}_{1}+\cdots +{\mathbb Z}_{p}\bar{a}_{n-x}$. Recall that $\langle a_{1},\ldots ,
a_{n-2x}\rangle G^{p}$ is a characteristic subgroup of $G$. Let $W={\mathbb Z}_{p}\bar{a}_{1}+\cdots +{\mathbb Z}_{p}\bar{a}_{n-2x}$ be the corresponding subspace
of $V$. Then let 
        $$H=\{\phi\in \mbox{GL}(n-x,p): \phi(W)=W\}.$$
There is now a natural action from $H$ on $Q(n,x)$. Suppose we have some presentation with generators $a_{1},\ldots ,a_{n-x}$ as above. Let $\phi\in H$ and suppose that 
$$ \overline{a_{i}}^{\phi}=\beta_{1}(i) \overline{a_{1}} + \cdots + \beta_{n-x}(i)\overline{a_{n-x}}.$$
We then get a new presentation in $Q(n,x)$ for $G$ with respect to the generators $b_{1},\ldots ,b_{n-x}$ where $b_{i}=a_{1}^{\beta_{1}(i)}\cdots a_{n-x}^{\beta_{n-x}(i)}$. \\ \\
Suppose there are $l$ powerfully nilpotent groups of exponent $p^{2}$ and order $p^{n}$ where $|G^{p}|=p^{x}$. Pick powerfully nilpotent presentations $p_{1},\ldots ,p_{l}\in {\mathcal P}(n,x)$ for these. Let $q$ be any powerfully
nilpotent presentation in ${\mathcal P}(n,x)$ of a group $K$ with generators $b_{1},\ldots ,b_{n-x}$. Then $q$ is a presentation for an isomorphic group $G$  with presentation $p_{i}$ and generators $a_{1},\ldots ,a_{n-x}$. Let $\phi:K \rightarrow G$ be an isomorphism and let $\psi:K/K^{p}
\rightarrow G/G^{p}$ be the corresponding linear automorphism. This gives us a linear automorphism $\tau\in H$ induced by $\tau(\bar{a}_{i})=\psi(\bar{b}_{i})$. Thus $q=p_{i}^{\tau}$. Thus the orbits of $p_{1},\ldots ,p_{l}$ cover ${\mathcal P}(n,x)$ and thus 
         $${\mathcal P}(n,x)\subseteq p_{1}^{H}\cup p_{2}^{H}\cup \cdots \cup p_{l}^{H}.$$
From this we get 
                             $$p^{\alpha n^{3}+o(n^{3})}=|{\mathcal P}(n,x)|\leq |p_{1}^{H}|+\cdots +|p_{l}^{H}|\leq l\cdot |H|\leq l\cdot p^{n^{2}},$$
and it follows that $l\geq p^{\alpha n^{3}+o(n^{3})}$. We thus get the following result. \\
\begin{theo} The number of powerfully nilpotent groups of 
exponent $p^{2}$ and order $p^{n}$ is 
           $p^{\alpha n^{3}+o(n^{3})},$
where $\alpha =\frac{9+4\sqrt{2}}{394}.$
\end{theo}
{\bf Remark}. Using a similar analysis one can estimate the growth of all powerful $p$-groups of exponent $p^{2}$ and order $p^{n}$, where $p$ is odd. This turns out to be
   $$p^{\frac{2}{27}n^{3}+o(n^{3})}.$$
Thus while the powerfully nilpotent $p$-groups of exponent $p^{2}$ and order $p^{n}$ are very numerous they are sparse within the larger class of all powerful $p$-groups of exponent $p^{2}$ and order $p^{n}$ as $n$ tends to infinity. \mbox{}\\ \\


\begin{thebibliography}{99}
%
\bibitem{l} A. Lubotzky and A. Mann, Powerful $p$-groups, {\it J. Algebra}, {\bf 105} (1987), 485--505.
%
\bibitem{m} P. Moravec and G. Traustason,
Powerful 2-Engel groups, {\it Comm. Algebra},  {\bf 36} (2008),
no. 11, 4096--4119.
%
\bibitem{sh} A. Shalev, On almost fixed point free automorphisms, {\it J. Algebra}, {\bf 157}, 271--282. 
%
\bibitem{s1} L. Sorkatti and G. Traustason,
Nilpotent symplectic alternating algebras {I}, {\it J. Algebra}  {\bf 423} (2014), 615--635.
%
%
\bibitem{tr1} G. Traustason,
Powerful 2-Engel groups II, {\it J. Algebra} {\bf 319} (2008), no. 8,
3301--3323.
%
\bibitem{tr2} G. Traustason,
Symplectic Alternating Algebras, {\it Internat. J. Algebra Comput.}
{\bf 18} (2008), no. 4, 719--757.
%
\bibitem{z}{GAP4}, The Gap~Group, \emph{GAP -- Groups, Algorithms, and Programming, Version 4.9.1}; 2018, \verb+(https://www.gap-system.org)+.
\end{thebibliography}
\end{document}